\documentclass[12pt]{article}
\usepackage{sfmath}

\newtheorem{theorem}{Theorem}
\newtheorem{lemma}[theorem]{Lemma}

\newtheorem{proposition}[theorem]{Proposition}

\newtheorem{definition}[theorem]{Definition}

\newtheorem{notation}[theorem]{Notation}
\newtheorem{remark}[theorem]{Remark}

\begin{document}

\title{On tropical and Kapranov ranks of tropical matrices}
\author{Elena Rubei}
\date{\hspace*{1cm}}
\maketitle

\vspace{-1.1cm}

{\small 






\bigskip

\section{Introduction}

Let  us consider the tropical semiring $({\bf R}, \oplus, \odot)$,
where:
$$ x \oplus y := min\{x,y\} \;\;\;\;\;\;\;\;  x \odot y := x+ y $$

For matrices with coefficients in the tropical semiring, we can define the 
tropical sum and the tropical product in the obvious way.

In \cite{DSS} Develin, Santos and Sturmfels compared three possible 
definitions for the rank of a tropical matrix.

\begin{definition} Let $A \in M(r \times r, {\bf R}) $. We define 
$$ det (A) := \oplus_{ \sigma \in S_r} a_{1, \sigma (1)} \odot....\odot
 a_{r, \sigma (r)} = min_{ \sigma \in S_r} 
\{a_{1, \sigma (1)} + .... + a_{r, \sigma (r)}\} $$
where $S_r$ denotes the symmetric group on $r$ elements.

Besides we say that $A$ is singular if the 
minimum in $det(A) $ is attained at least twice.
\end{definition}

\begin{definition} Let $M$ be  a matrix $m\times n$ with entries in
 ${\bf R}$. We define:

1) (Barvinok rank) $rk_B(M) := min\{ r | \;  \exists A_1, ..., A_r \in M(
  m \times 1, {\bf R}) \; and \;  B_1, ..., B_r \in M (1 \times n, {\bf R}) 
\; s.t. \; M =   A_1 \odot B_1 \oplus ...\oplus A_r \odot B_r\}$

2) (Kapranov rank)  $rk_K(M) := min\{ r | \; \exists \;a\; tropical \;
 linear\; space \; of \; dimension \; r 
$ $  containing $ $  the$ $ columns\; of\; M\}$ (see \S2 for the definition
 of tropical linear space)

3) (Tropical rank) $rk_t (M) := max \{r |\; \exists \; a \;nonsingular\; 
 minor \; r \times r \; of \; M \} $ 

\end{definition}

Develin, Santos and Sturmfels proved

\begin{theorem} \label{dis}
 \cite{DSS}. For every matrix $M \in M( m \times n, {\bf R})$ we have
$$ 1 \leq rk_{t} (M) \leq rk_{K} (M) \leq rk_{B} (M) \leq  min\{m,n\} $$
Besides, if either $rk_{t} (M)$ or $rk_{K} (M)$ 
 is equal either to $1$ or to $ 2$ or to $ min\{m,n\}$,
then $rk_{K} (M) = rk_{t} (M)$.
\end{theorem}

Furthermore they showed that the inequalities may be strict: to show that
 the first inequality may be strict, they exhibit a matrix $7 \times 7$ 
whose tropical  rank is $3$ and whose Kapranov rank is $4$; they 
wonder whether there exists a matrix $5 \times 5$ whose tropical  and 
Kapranov
 ranks disagree, in particular with tropical rank $3$ ad Kapranov rank
strictly greater than $3$ (observe that if $M$ is a matrix $5 \times 5$ 
with tropical rank $1,2,4,5$, then, by the theorem above and Remark
\ref{k+1k+1} below, which is an obvious corollary of the theorem above,  
Kapranov rank is equal
to the tropical rank; 
then the unique case for a matrix $5 \times 5$ when 
 tropical rank and Kapranov rank may disagree is when the tropical 
rank is $3$).

Here we prove the following theorem,

\begin{theorem} \label{mio} 
Let $A$ be a matrix $g \times 5$ of tropical rank $3$. 
Then also the Kapranov rank is  $3$. 
\end{theorem}

Finally we quote the following result  by Kim and Roush:

\begin{theorem} \label{KR}
\cite{KR}. Kapranov rank of tropical matrices is not bounded in terms
of tropical ranks, that is, given a positive integer $r$, it is impossible 
to find a bound $N(r)$ so that all matrices of tropical rank $r$ 
have Kapranov rank less or equal than $N(r)$. 
\end{theorem}

\section{Some recalls}

Before proving our theorem, we need to recall what a linear tropical 
space is 
(we have used it in the definition of Kapranov rank) and a proposition 
by Develin, Santos and Sturmfels.

Let $K$ be the algebraic closure of $ {\bf C} (\tau)$, that is the field of
 Puiseux  series with complex coefficients.
\begin{eqnarray*} 
 K = \{c_1 \tau^{a_1} + c_2 
\tau^{a_2}+...... | \; c_i \in {\bf C}, \; a_i \in 
{\bf Q},\; a_1 < a_2  <....\;and\; the \;  a_i \; have\; \\
 a\; common \; denominator \}
\end{eqnarray*}
Let us define
$$ \tilde{K} = 
\{c_1 \tau^{a_1} + c_2 \tau^{a_2}+...... | \; c_i \in {\bf C}, \; a_i \in 
{\bf R},\; a_1 < a_2  <....\}$$
Let $ord : \tilde{K} \rightarrow {\bf R}$
be the map 
$$  (c_1 \tau^{a_1} + c_2 \tau^{a_2}+......) \mapsto a_1 $$
and we call $ord $ also the obvious map 
$ord : \tilde{K}^d \rightarrow {\bf R}^d$

Let $orc : \tilde{K} \rightarrow {\bf R}$
be the map 
$$  (c_1 \tau^{a_1} + c_2 \tau^{a_2}+......) \mapsto c_1 $$

Let $I$ be an ideal in $K [x_1,..., x_d]$.
Let $$ \tilde{V} (I) = \{x \in (\tilde{K}-\{0\})^d | \; f(x) =0 \; 
\forall f \in I \} $$ 
We define the tropical variety of $I$ as $ord (\tilde{V}(I))$.

We say that a tropical variety $ord(\tilde{V}(I) ) $ is linear if $I$ has a 
set of linear generators.

\begin{remark} \label{k+1k+1} Let $M$ be a matrix $(k+1) \times n $ 
with $ rk_t (M) =k$ with $ n \geq k+1$. Then $ rk_K (M) = k$.
\end{remark}

{\it Proof.} 
If $rk_K (M) $ were not $k$, it would be $ k+1$ (because $rk_K (M) \geq 
rk_t (M) =k$ by the first disequality of
Theorem \ref{dis}), but in this case also $rk_t (M) $ would be
$k+1$ by the second statement of Theorem \ref{dis}.
\hfill \framebox(7,7)

\begin{proposition} \label{lift} 
\cite{DSS} Let $M \in M( m \times n, {\bf R})$.
The following are equivalent:

1) $rk_K (M) \leq r$

2) there exists $F \in M(m \times n, \tilde{K}- \{0\}) $ 
with $rk (F) \leq r$ and $F$ ``lift'' of $M$, that is $ ord (F) =M$.

\end{proposition}

\section{Proof of the theorem}

\begin{notation}
$\bullet $ 
When we write a matrix, if in an entry we write $\ast$, this 
means that this 
entry is nonzero.  

$\bullet$ U.s.r.c. means ``up to swapping rows and columns''.

$\bullet$ If $A$ is a matrix, $A_{\hat{\i},\cdot} $ is the 
matrix obtained from 
$A$ by taking off the $i$-th row, $A_{(i)}$ is the $i$-th row 
and  $A^{(j)}$ is the $j$-th column.

$\bullet$ We say that we can develop a tropical matrix $A$ $m \times n$ 
by one of 
its rows $A_{(i)}$ with coefficients of order 
$ l_1,..., l_{m-1}$ if, for any 
lift $F$ of $A_{\hat{\i},\cdot} $, we can find coefficients 
in $K$ of order 
$ l_1,..., l_{m-1}$ such that the linear combination of the rows of $F$ 
with these coefficients is a lift of $A_{(i)}$.

\end{notation}

\begin{remark} \label{IR} {\bf (Important remark).}
Suppose  we know that every matrix $(g-1) \times n$ with tropical 
rank $s$ has the Kapranov rank $s$ for any $s \leq k$.
Let $A$ be a matrix $g \times n$  with tropical rank  $k$ such that
we can develop $A$ or one of its submatrix   $p 
\times n$ by one of its row $A_{(i)}$.
Then also Kapranov rank of $A$ is $k$.   

In fact there exists a lift of rank $ \leq k$  of $A_{\hat{\i},\cdot}$;
this together with the linear 
combination of its rows given by the developing will give a lift  
of rank $ \leq k$ of $A$. 

We will apply this remark with $n=5$.
\end{remark}

\begin{remark} \label{lemmino}
Let $a,b,c,d \in K$ of order $0$ and $ r, r' \in {\bf R}^{+}$ with 
$r \neq r'$. Then there exist $ \lambda_1, \lambda_2 \in K$ of order greater 
or  equal than $-m$ ($m \in {\bf R}^{+}$) such that 
$$ord ( \lambda_1 a + \lambda_2  b)=r \;\;\;\;\;\;\;\;
ord ( \lambda_1 c + \lambda_2  d)=r' $$ 
iff 
{\small $\left(
\begin{array}{c}
a  \\ b
\end{array}
\right)$}
and  {\small $\left(\begin{array}{c}
c  \\ d
\end{array}
\right)$}
are indipendent and $ord (a d -b c ) \leq min\{r,r'\}+m$
\end{remark}

\begin{remark} \label{ande}
Let $ r \in {\bf N}$, $r \geq 2$. Let $ X_1,..., X_r \in K$ (and possibly 
$ Z_1,..., Z_s \in K$)  
of order respectively $ x_1,...., x_r $,  ($z_1,..., z_s$).  
Let $l_1,..., l_r, k,(m_1,..., m_s)$ $\in K$ such that 
$$ h=:x_1+l_1= ...=x_r+l_r \leq k\;\;\; \; \; (m_1+z_1 ,..., m_s+z_s > h)$$  
a) 
Then there exist $ \lambda_1,..., \lambda_r, (\mu_1,..., \mu_s)$ of order 
$ l_1,..., l_r, (m_1,..., m_s)$ such that 
$$ \lambda_1 X_1 +...... \lambda_r X_r (+ \mu_1 Z_1 +....+ \mu_s Z_s)$$ 
has order $k$.

b) If $ k> h$, choosen $\lambda_1,...., \lambda_{r-1} $ of order
$ l_1,...., l_{r-1}$ such that $ \lambda_1 X_1 + .... \lambda_{r-1} 
X_{r-1}$ has order $h$ (thus if $r=2$ for any choice of $\lambda_1 $ of 
order $l_1$) 
(and choosen $ \mu_1,..., \mu_s$ of order $m_1,..., m_s$), we can choose 
$\lambda_r$ of order $ l_r$ such that 
$$ \lambda_1 X_1 + .... \lambda_r X_r (+ \mu_1 Z_1 +....+ \mu_s Z_s)$$ 
has order  $k$ 

c) If $k= h$, choosen in any way 
 $\lambda_1,...., \lambda_{r-1} $ of order
$ l_1,...., l_{r-1}$,
(and choosen $ \mu_1,..., \mu_s$ of order $m_1,..., m_s$), we can choose 
$\lambda_r$ of order $ l_r$ such that 
$$ \lambda_1 X_1 + .... \lambda_r X_r (+ \mu_1 Z_1 +....+ \mu_s Z_s)$$ 
has order  $k$.

\end{remark}

\begin{remark}\label{remquad}
Let $Y_1,Y_2, Y_3, B_1, B_2, Z \in K$ of order respectively
$y_1, y_2, y_3 , b, b, z$. Let $ t \geq z > b$. Then there exist 
$  \lambda_1, \lambda_2 \lambda_3 \in K $ of order $0$ 
such that 

$  \lambda_1 Y_1 + \lambda_2 Y_2 + \lambda_3 Y_3$ has order $y$ 

$  \lambda_1 B_1 + \lambda_2  B_2 + \lambda_3 Z $ has order $t$. 

\end{remark}

{\it Proof.}
Obviously if $ord (\lambda_i)= 0$ $i=1,2,3$,  then 
$  \lambda_1 Y_1 + \lambda_2 Y_2 + \lambda_3 Y_3$ has order $y$ iff 
$ \sum_{i=1,2,3} orc (\lambda_i ) orc( Y_i) \neq 0$. 

Besides observe that if $ord (\lambda_i)= 0$ $i=1,2,3$ and 
if  $  \lambda_1 B_1 + \lambda_2  B_2 + 
\lambda_3 Z $ has order $t$ then 
$ orc ( \lambda_2 ) = orc (\lambda_1) \frac{orc(B_1)}{orc (B_2)} $.

So choose $ \lambda_1$ and $ \lambda_3 $ of ord $0$ such that 
$ orc (\lambda_1 ) orc( Y_1)  + 
orc ( \lambda_1 )  \frac{orc(B_1)}{orc (B_2)} + 
orc ( \lambda_3 ) orc (Y_3) \neq 0$. 
Then, by using the equation $ ord ( \lambda_1 B_1 + \lambda_2  B_2 + 
\lambda_3 Z) = t $, define $\lambda_2 $, for instance take $ \lambda_2
=  \frac{-\lambda_1 B_1  -\lambda_3 Z + \tau^t }{B_2}$
(it will  have order $0$). 
 
\hfill \framebox(7,7)

\begin{remark}\label{remquadsbar}
Let $A_1,A_2, B_2, B_3, H, Y \in K$ of order respectively
$a,a, b, b ,h, y$. Let $ t , h > b$ and $ z , y > a$. Then there 
exist 
$  \lambda_1, \lambda_2 \lambda_3 \in K $ of order $0$ 
such that 

$  \lambda_1 A_1 + \lambda_2 A_2 + \lambda_3 Y $ has order $z$ 

$  \lambda_1 H + \lambda_2  B_2 + \lambda_3 B_3 $ has order $t$. 

\end{remark}

{\it Proof.}
Choose $\lambda_1$ and $ \lambda_3$ such that 
$$ \frac{-\lambda_1 A_1 - \lambda_3 Y + \tau^z }{A_2}= 
\frac{-\lambda_1 H - \lambda_3 B_3 + \tau^t }{B_2}$$
Then define $ \lambda_2 $ as  one of the two members of the above equation.
\hfill \framebox(7,7)

\begin{remark}\label{sbis}
Let $A_1,A_2, B_2, B_3, B_4, H, Y,W \in K$ of order respectively
$a,a, b, b, b ,h, y,w$. Let $ t , h > b$ and $ z , y ,w> a$. 
Let $ c_1,c_3 c_4 ,c \in {\bf R} -\{0\}$. Then there 
exist 
$  \lambda_1, \lambda_2 \lambda_3 , \lambda_4 \in K $ of order $0$ 
such that 

$  \lambda_1 A_1 + \lambda_2 A_2 + \lambda_3 Y + \lambda_4 W$ has order $z$ 

$  \lambda_1 H + \lambda_2  B_2 + \lambda_3 B_3  + \lambda_4 B_4$ 
has order $t$ 

$c_1 orc(\lambda_1) + c_3 orc(\lambda_3) +
c_4 orc(\lambda_4) \neq 0$ 
 
$ orc(\lambda_1)/ orc(\lambda_4) \neq c$.

\end{remark}

\begin{lemma} \label{5linedup}
Let $A$ be a  a matrix $5 \times 5$ of tropical rank $3$. Suppose that there
is a row with $ 5$ zeroes, that every entry is greater or equal than 
$0$,   that in every column there are at least two zeroes. 
Then also Kapranov rank of $A$  is  $3$. 
\end{lemma}

{\it Proof.}
Obviously we can develop by the zero row.
\hfill \framebox(7,7)

\begin{lemma}\label{caso3322}
Let $A$ be a matrix $4 \times 5$ of kind 
$$\left(\begin{array}{ccccc}
0 & 0 & 0 & . & . \\
. & . & . & 0 & . \\
. & . & . & 0 & 0 \\
. & . & . & 0 & 0 
\end{array}
\right)$$
with all coefficients greater or equal than $0$. 
If the tropical rank is $3$, then we can develop 
$A$  by one of its row.
\end{lemma}

{\it Proof.}
Observe that the tropical determinant of a matrix 
$$\left(\begin{array}{cccc}
0 & 0 & 0 & .  \\
. & . & . & 0  \\
. & . & . & 0  \\
. & . & . & 0 
\end{array}
\right)$$
with all coefficient greater or equal than $0$, is the minimum $M$ 
of the sums of the couples of 
coefficients  of  the left lower submatrix $3 \times 3$ which 
are not in the same row and column. 
So, if the tropical rank is less or equal than $3$, there 
are at least two couples of coefficients, not lined up, 
whose sum is $M$. 
There are  two possibilities:

1) these two couples are disjoint, 
$a,b+h$ and  $a+h,b$ 
let's say $h \in {\bf R}^+$, $a \leq b$.

2) these two couples have an element in common, that is there are two 
coefficients equal to $b$ and one $a$ not lined up, whose sum is $M$.

1) There are 3 subcases, u.s.r.c. and up to transposing,

$$ 1.1) \left(\begin{array}{cccc}
0   & 0 & 0 & .  \\
a+h & a & . & 0  \\
.   & b & . & 0  \\
b+h & . & . & 0  
\end{array}
\right) \;\;\;\; 1.2) \left(\begin{array}{cccc}
0 & 0 & 0 & .  \\
a+h & a & . & 0 \\
. & b & . &  0  \\
. & . & b+h & 0  
\end{array}
\right) \;\;\;\; 1.3)  \left(\begin{array}{cccc}
0 & 0 & 0 & .  \\
a+h & a & . & 0  \\
b+h & b & . & 0  \\
. & . & . & 0  
\end{array}
\right)$$

2) U.s.r.c. and up to transposing, we have two cases:

$$ 2.1) \left(\begin{array}{cccc}
0 & 0 & 0 & .  \\
. & a & . & 0  \\
. & . & b & 0  \\
. & . & b & 0  
\end{array}
\right) \;\;\;\;\; 2.2) \left(\begin{array}{cccc}
0 & 0 & 0 & .   \\
b & . & . & 0   \\
. & b & . & 0   \\
. & . & a & 0  
\end{array}
\right)$$
with $a \leq b$

Observe that  the submatrix $A_{\cdot, (1,2,3,4)}$ of $A$ must be of 
one of the forms above, u.s.r.c., since its rank is less or equal than $3$.

Case 1.1) Since also $A_{\cdot , (1,2,3,5)}$ 
has rank less or equal than $3$
then, u.s.c., $$A=\left(\begin{array}{ccccc}
0 & 0 & 0 & . & . \\
a+h & a & u & 0 & . \\
s & b & y & 0 & 0 \\
b+h & t & w & 0 & 0 
\end{array}
\right)\;\;\; or \;\;\;A=\left(\begin{array}{ccccc}
0 & 0 & 0 & . & . \\
a+h & a & u & 0 & . \\
b+h & t & w & 0 & 0 \\
s & b & y & 0 & 0 
\end{array}
\right)$$
(that is $a+h$ and $a$ must be lined up with the nonzero element 
among $ A_{2,5},  A_{3,5}, A_{4,5}$). U.s.r., we can suppose the first.

We must have $t \geq b$, $ w,s,y \geq b+h$, $ u \geq a+h$.
Suppose the inequalities are strict.

we can divide into many cases, we use Remark \ref{ande}:

\medskip

$\bullet $ $ t < min\{w,y, u+b-a\}$ 

- if $ y \neq  u+b-a$, we can develop
$A$ by $A_{(4)}$ with 
coefficients   of order  respectively $ min\{w,y, u+b-a\}, b-a, 0$

- if $ y =  u+b-a$:

if $ w \leq y$ 
we can develop
$A$ by $A_{(4)}$ with 
coefficients   of order  respectively $ w, b-a, 0$;

if $ w > y$ 
we can develop
$A_{(2,3,4) , \cdot}$ by $A_{(3)}$ with 
coefficients   of order  respectively $  b-a, 0$;

\medskip

$\bullet$  $ t \geq min\{w,y, u+b-a\}$

- if  $ y \neq  u+b-a$,
  we can develop $A$ 
by $A_{(4)}$ with 
coefficients   of order  respectively $ min\{w,y, u+b-a\}, b-a, 0$ (use 
Remark \ref{remquadsbar});

- if  $ y =  u+b-a$,

if $w \leq y$ we can develop
$A$ by $A_{(4)}$ with 
coefficients   of order  respectively $ w, b-a, 0$ (if $w=y$ use Remark 
\ref{remquad});

if $w>y$ 
we can develop
$A_{(2,3,4) , \cdot}$ by $A_{(4)}$ with 
coefficients   of order  respectively $ b-a, 0$;



Analogously the other cases. 
\hfill \framebox(7,7)

\begin{lemma}\label{casospecchio}
Let $A$ be a matrix $g \times 5$ of kind 
$$\left(\begin{array}{ccccc}
v' & u' & r' & 0 & 0 \\
v & u & r & 0 & 0 \\
0 & 0 & 0 & . & . \\
\vdots & \vdots & \vdots & \vdots & \vdots \\
0 & 0 & 0 & . & . 
\end{array}
\right)$$
with all coefficients greater or equal than $0$.
If  the tropical rank is $3$, then also Kapranov rank is  $3$. 
\end{lemma}

{\it Proof.}
We can suppose that $ r = min \{v,v', u,u',r',r\}$.

Choose  lifts $  F^{(1)}, F^{(4)}, F^{(5)}$   
of $ A^{(1)}, A^{(4)}, A^{(5)}$ 
of the kind respectively:
$$ \left(\begin{array}{c}
c \tau^{v'} \\
h \tau^{v} \\
l_1 \\
 \vdots \\
l_{g-2}
\end{array}
\right) \;\;\; 
\left(\begin{array}{c}
1 \\
1 \\
\cdot \\
\vdots \\
\cdot
\end{array}
\right) \;\;\;\left(\begin{array}{c}
-1 + 3t^{min \{v',r'\}} \\
-1 + 2 t^{r} \\
\cdot \\
\vdots \\
\cdot
\end{array}
\right)  $$
with $c \in K$ of order $0$ and $ h,l_1,...,l_{g-2} \in {\bf R} - \{0\}$

By Remark \ref{lemmino} 
 there exist $  \mu_2, \mu_3 \in K$ of order $ \geq 0$ 
such that 
$$ ord \left(  \mu_2 
\left( \begin{array}{c} 
1\\
1
\end{array} \right)
+ \mu_3 
\left( \begin{array}{c} 
-1 + 3t^{min \{v',r'\}} \\
-1 + 2 t^{r} 
\end{array} \right)
\right)= 
\left( \begin{array}{c} 
min\{v',u'\} \\
u
\end{array} \right)$$

Choose $ c $ such that  $  (F^{(1)}+ F^{(4)}+ F^{(5)})_1$ has $ord$ 
equal to  
$ A^{(3)}_1$ that is $r'$.

There exist $\mu_1 \in K $ of order $0$ such that  $(\mu_1 
F^{(1)}+ \mu_2 F^{(4)}+ \mu_3 F^{(5)})_1$ has $ ord$ equal to  
$ A^{(2)}_1$, that is $u'$.

Obviously we can choose $ h,l_1,...,l_{g-2}$ such that 
$$  ord (F^{(1)}+ F^{(4)}+ F^{(5)}) = A^{(3)} \;\;\;\;and \;\;\;\; 
 ord (\mu_1 F^{(1)}+ \mu_2 F^{(4)}+ \mu_3 F^{(5)}) = A^{(2)}$$
\hfill \framebox(7,7)

\begin{lemma} \label{4linedup}
Let $A$ be a  a matrix $5 \times 5$ of tropical rank $3$. Suppose that there
is a row with $4$ zeroes, that every entry is greater or equal than 
$0$,   that in every column there are at least two zeroes and that 
in every row there is at least one zero. 
Then, if we are not in case of Lemma \ref{casospecchio}, we can develop
$A$ by the row with $4$ zeroes
and then also the Kapranov rank of $A$  is  $3$. 
\end{lemma}

{\bf Proof.} U.s.r.c. we can suppose that $A_{1,2}= ....= A_{1,5}=0 $.  
U.s.r, we can suppose $A_{2,1}= ....= A_{k,1}=0 $,  
$0 < A_{k+1,1}, ...., A_{s,1} \leq r $,
 $A_{s+1,1}, ...., A_{5,1} >r $ with $k \geq 3$.

We state that we can develop $A$ by the first row with coefficients
 $ \lambda_2,..., 
\lambda_5$ of order $0$, in fact:

 let $(F_{i,j})_{i=2,...,4\;\;\; j=1,..., 5} $ be a 
lift of $ A_{[2,...,5] , \cdot }$

i) to have $ ord ( \sum_{i=2,...,5} \lambda_i F_{i,1} ) =r $, it 
is sufficient to choose $\lambda_2,..., \lambda_5$ of order $0$ such that 
$$ \sum_{i=3,..., k } orc (\lambda_i ) orc( F_{i,1}) \neq 0 $$
  (so that 
$ ord( \sum_{i=3,..., k } \lambda_i F_{i,1})= 0$)
and then choose $$ \lambda_2 = \frac{ -\sum_{i=3,..., s } 
\lambda_i F_{i,1} +\tau^r}{F_{2,1}}$$
ii)  to have $ ord ( \sum_{i=2,...,5} \lambda_i F_{i,j} ) =0 $ 
for some $ j \in \{2,...,5\}$, it is sufficient
  to choose $\lambda_2,..., \lambda_5$ of order $0$ such that 
$ \sum_{i \in \{2,...,5\} \;s.t.\;A_{i,j} = 0 } 
orc (\lambda_i ) orc( F_{i,j}) \neq 0 $ 

Obviously we can easily 
find  $\lambda_2,..., \lambda_5$ satisfying either 
i) and ii) unless we are in one of the following cases (u.s.r.c.)
$$A=\left(\begin{array}{ccccc}
r & 0 & 0 & 0 & 0 \\
0 & 0 & . & . & . \\
0 & 0 & . & . & . \\
0 & 0 & . & . & . \\
0 & 0 & . & . & . 
\end{array}
\right) \;\;\; A=\left(\begin{array}{ccccc}
r & 0 & 0 & 0 & 0 \\
0 & 0 & . & . & . \\
0 & 0 & . & . & . \\
0 & 0 & . & . & . \\
\ast & \ast & . & . & . 
\end{array}
\right) \;\; \; A=\left(\begin{array}{ccccc}
r & 0 & 0 & 0 & 0 \\
0 & 0 & . & . & . \\
0 & 0 & . & . & . \\
\ast & \ast & . & . & . \\
\ast & \ast & . & . & . 
\end{array}
\right)$$
In the first two cases we can conclude by Lemma \ref{caso3322}. 
In the last one
observe that, by assumption, we must have a zero in each of the last 
two rows. Observe that 
if $$ A=\left(\begin{array}{ccccc}
r & 0 & 0 & 0 & 0 \\
0 & 0 & . & . & . \\
0 & 0 & . & . & . \\
\ast & \ast & 0 & . & . \\
\ast & \ast & \ast & 0 & . 
\end{array}
\right)$$
we would have that the tropical rank of $ A_{\hat{3}, \hat{5}}$ is $4$, 
 which contradicts the assumption $rk_t (A) =3$.
Thus the only possibilities are (u.s.c.)
$$ A=\left(\begin{array}{ccccc}
r & 0 & 0 & 0 & 0 \\
0 & 0 & . & . & . \\
0 & 0 & . & . & . \\
\ast & \ast & 0 & \ast & \ast \\
\ast & \ast & 0 & \ast & \ast 
\end{array}
\right), \;\; \;A=\left(\begin{array}{ccccc}
r & 0 & 0 & 0 & 0 \\
0 & 0 & . & . & . \\
0 & 0 & . & . & . \\
\ast & \ast & 0 & 0 & \ast \\
\ast & \ast & 0 & 0 & \ast 
\end{array}
\right)\;\;\; A=\left(\begin{array}{ccccc}
r & 0 & 0 & 0 & 0 \\
0 & 0 & . & . & . \\
0 & 0 & . & . & . \\
\ast & \ast & 0 & 0 & 0 \\
\ast & \ast & 0 & 0 &  0
\end{array}
\right) $$
In the last case we are in case of 
Lemma \ref{casospecchio}. In the first 
two observe that at least one of $ A_{2,5}$ and $ A_{3,5}$ must 
be zero, so we can conclude by Lemma \ref{caso3322}.
\hfill \framebox(7,7)

\bigskip

{\bf Proof of Theorem \ref{mio}.}
Induction on $g$. The case $g=4$ follows from Remark \ref{k+1k+1}. So
 suppose to know the statement for $g-1$. Let $A$ be $ g \times 5$. We will
show that in almost all cases we can develop $A$ or one of its 
submatrix $ p \times 5$ by one of its rows and we will conclude by Remark 
\ref{IR}.

U.s.r. we can suppose that the submatrix of $A$ given by the first 
$4$ rows,
$A_{[1,..,4], \cdot} $ has tropical rank $3$.

Since $ rk_t ( A_{[1,..,4], \cdot}) =3$ 
then also $ rk_K ( A_{[1,..,4], \cdot}) =3$ 
(in fact if Kapranov rank were $4$, then also the tropical rank would be $4$
by Theorem \ref{dis}).

So we can suppose that the columns of 
$A_{[1,..,4], \cdot} $ are in a tropical hyperplane. 
We can suppose that this 
hyperplane is given by $ x_1 \oplus ... \oplus x_4$.
So, up to summing some constants to the columns of $A$, we can suppose that 
in all the columns of $A_{[1,..,4], \cdot} $ we have two zero entries and 
all the other entries are greater or equal than $0$. 

If two zero entries are in the same column we say that they are twin zeroes.

For any $i=5,...,g$, 
up to summing a constant to all the entries of $A_{(i)}$,
we can suppose also that the entries of $A_{(i)}$ are greater or equal 
than $0$ and that the minimum is $0$.

We can distinguish in five cases according to the maximum number $m$ of the 
couples of twin zeroes in the same two rows of $A_{[1,2,3,4], \cdot}$.

$\bullet$ $m=1$.

Observe that in this 
case there exists a row among the first four rows of $A$ 
with at least $3$ zeroes
(in fact if every row of $A_{[1,..,4], \cdot} $ had at most $2$ zeroes, then 
the number of the zeroes of $A_{[1,..,4], \cdot} $ 
would be less or equal than $8$, but the number 
of the zeroes  of $A_{[1,..,4], \cdot} $ is at least $2 \cdot 5=10$ since 
every columns contains at least $2$ zeroes).
So, u.s.r.c, we can suppose
$A_{1,1}= A_{1,2}=A_{1,3}=0$.

Since $m=1$,
$A$ must be like this (u.s.r.c.):
$$\left(\begin{array}{ccccc}
0 & 0 & 0       & a & b \\
0 & \ast & \ast & 0 & b' \\
\ast & 0 & \ast & 0 & 0 \\
\ast & \ast & 0 & a' & 0 \\
\vdots & \vdots & \vdots & \vdots & \vdots
\end{array}
\right)$$
with $a,b,a',b' \neq 0$.
We can develop $A_{[1,2,3,4], \cdot}$ by the first row with coefficients
of order $0$: if either $ a' > a$ or $ b'> b$, use 
Remark \ref{ande}, otherwise  use Remark \ref{remquadsbar}.

$\bullet$ $m=5$. It follows from Lemma \ref{5linedup}.

$\bullet$ $m=4$. It follows from Lemma \ref{5linedup} or Lemma 
 \ref{4linedup} or Lemma \ref{casospecchio}.

$\bullet$ $m=3$. 
U.s.r.c.,  we can suppose 
that  $A_{i,j}=0$ for $ i=1,2$, $j=1,2,3$. 

By Lemmas \ref{5linedup}, \ref{4linedup} and \ref{casospecchio}, 
we can suppose that
 $A$ is like this:
$$\left(\begin{array}{ccccc}
0 & 0 & 0 & v & v' \\
0 & 0 & 0 & r & r' \\
. & . & . & . & . \\
. & . & . & . & . \\
\vdots & \vdots & \vdots & \vdots & \vdots
\end{array}
\right)$$
with $r,r',v,v' \neq 0$ and, u.s.r.c., we can suppose 
$r= min \{r,r',v,v'\}$.

We can distinguish in two main  cases: 

1) the minimum in $\{r,r',v,v'\}$ is attained at least twice.

2) the minimum in $\{r,r',v,v'\}$ is attained only once.

\bigskip

1) The minimum in $\{r,r',v,v'\}$ is attained at least twice.

Obviously we can suppose that $(r,r') \neq (v,v')$.
If 
$$A=\left(\begin{array}{ccccc}
0 & 0 & 0 & v & v'  \\
0 & 0 & 0 & r & r \\
. & . & . & 0 & 0 \\
. & . & . & 0 & 0 \\
\vdots & \vdots & \vdots & \vdots & \vdots
\end{array}
\right) \;\; \;  or \; \;\; A=\left(\begin{array}{ccccc}
0 & 0 & 0 & v & r  \\
0 & 0 & 0 & r & r' \\
. & . & . & 0 & 0 \\
. & . & . & 0 & 0 \\
\vdots & \vdots & \vdots & \vdots & \vdots
\end{array}
\right)$$ 
we can develop $A_{[1,2,3],\cdot }$ by $A_{(2)}$ with coefficients
of order $ 0$ for the first row and $ r$ for the third row
(in the second case we can suppose $v \neq r$, if not we are in the previous 
case).

If 
$$A=\left(\begin{array}{ccccc}
0 & 0 & 0 & r &v'  \\
0 & 0 & 0 & r & r' \\
. & . & . & 0 & 0 \\
. & . & . & 0 & 0 \\
\vdots & \vdots & \vdots & \vdots & \vdots
\end{array}
\right)$$ with $r \neq r'$,  we can suppose, u.s.r.c., that $v'\geq r'$;  
develop $A_{[1,2,3],\cdot }$ by   $A_{(2)}$ with coefficients
of order $ 0$ for the first row and $ r'$ for the third row.

2)  The minimum in $\{r,r',v,v'\}$ is attained only once.
Then the tropical rank of $\left(\begin{array}{ccc}
0 & v & v'  \\
0 & r & r' \\
. & 0 & 0 
\end{array}
\right)$ is $3$; besides observe that 
if we take off from $A$ one of the first three columns, then in every of the 
first four rows of the so obtained matrix
there are still at least two zeroes. So its rows 
 must stay in a unique hypersurface and this hypersurface must be $ x_1
\oplus x_2 \oplus x_3 \oplus x_4$. Therefore, for every $ i=5,...,g$,  
in $A_{(i)}$ there are at least two zeroes in the entries
$ 2,3,4,5$, two zeroes in the entries $ 1,2,4,5$ and two zeroes
in the entries $ 1,3,4,5$; thus, for every $ i=5,...,g$, either  
 $A_{(i)}=(.,.,.,0,0) $ or in  $A_{(i)}$ there at least two 
zeroes in the first three entries and one in the last two
or $ A_{(i)}=(0,0,0,.,.) $.
If there exists $ i \in \{5,...,g\}$ such that 
either $A_{(i)}$ is $ (.,.,.,0,0) $
or in  $A_{(i)}$
 there at least two  zeroes in the first three entries and one in the last two,
we can conclude by Lemma \ref{caso3322}. If, for any $ i=5,...,g$,
 $A_{(i)} =(0,0,0,.,.)  $
we can conclude by Lemma \ref{casospecchio}).

$\bullet$ $m=2$.

U.s.r.c.,  we can suppose 
that  $A_{i,j}= 0$ for $i,j \in \{1,2\}$.

We can divide into many cases according how many zeroes are among the 
six entries in the first two  lines and last three columns.

If  among the six entries
in the first two lines and last three columns, there were no zeroes, then  
 $A$ would be  like this:
$$\left(\begin{array}{ccccc}
0 & 0 & \ast & \ast & \ast \\
0 & 0 & \ast & \ast & \ast \\
. & . & 0 & 0 & 0 \\
. & . & 0 & 0 & 0 \\
\vdots & \vdots & \vdots & \vdots & \vdots
\end{array}
\right)$$
and we could conclude since $m$ should be $3$.

If  among the six entries
in the first two lines and last three columns, there are at least  three 
zeroes, then at least two of them must be lined up and so we can conclude by 
Lemma \ref{4linedup} and Lemma \ref{casospecchio}.

So we can suppose that  among the six entries
in the first two lines and last three columns,
there are one or two zeroes and, in the case of two zeroes, 
that they are not lined up.

So we can distinguish the following cases (u.s.r.c.):
$$ a) \;\; A= \left(\begin{array}{ccccc}
0 & 0 & \ast & \ast & \ast \\
0 & 0 & \ast & \ast & 0 \\
. & . & 0 & 0 & 0 \\
. & . & 0 & 0 & . \\
\vdots & \vdots & \vdots & \vdots & \vdots
\end{array}
\right) \;\;\;\; b)\; \; A= \left(\begin{array}{ccccc}
0 & 0 & 0 & \ast & \ast \\
0 & 0 & \ast & 0 & \ast \\
. & . & 0 & . & 0 \\
. & . & . & 0 & 0 \\
\vdots & \vdots & \vdots & \vdots & \vdots
\end{array}
\right)  \;\;\;\;
c)\; \; A= \left(\begin{array}{ccccc}
0 & 0 & 0 & \ast & \ast \\
0 & 0 & \ast & 0 & \ast \\
. & . & 0 & 0 & 0 \\
. & . & . & . & 0 \\
\vdots & \vdots & \vdots & \vdots & \vdots
\end{array}
\right)$$

\underline{Case a)}  $A= \left(\begin{array}{ccccc}
0 & 0 & \ast & \ast & \ast \\
0 & 0 & \ast & \ast & 0 \\
. & . & 0 & 0 & 0 \\
. & . & 0 & 0 & . \\
\vdots & \vdots & \vdots & \vdots & \vdots 
\end{array}
\right)$

We can suppose  $A= \left(\begin{array}{ccccc}
0 & 0 & \ast & \ast & \ast \\
0 & 0 & \ast & \ast & 0 \\
\ast & \ast & 0 & 0 & 0 \\
\ast & \ast & 0 & 0 & . \\
\vdots & \vdots & \vdots & \vdots & \vdots 
\end{array}
\right)$ by Lemma \ref{4linedup} and Lemma \ref{caso3322}.

For any $i=5,..., g$, in  $A_{(i)}$ there is at least a zero.
If a zero is in the first four entries, then we can conclude 
by Lemma \ref{caso3322}.
So the unique problem is when for every $ i=5,.., g$, 
 there is only one  zero 
in  $A_{(i)}$ and it is  the last entry.

Let 
$S_i $ be the submatrix   $A_{[1,2,i], [3,4]}$.
The tropical determinant of $A_{[124i], [2345]}
=  \left(\begin{array}{cccc}
0 & \ast & \ast & \ast \\
0 & \ast  & \ast & 0 \\
. & 0  & 0 & . \\
\ast & \ast  & \ast & 0
\end{array}
\right)$ is the 
minimum in $S_i$.
So  this minumum must be attained twice.

\medskip

If the minimum in $S_i$ is attained twice in the same column, 
say in the first,  we can develop
$A_{[1,2,3,i], \cdot}$  by the first or second row with coeffcients 
$ 0,min \{A_{1,4} , A_{2,4}, A_{i,4} \},0 $; while 
if the minimum $s_i$ in $S_i$ is attained twice but not in the same column
and not in the same row, 
  we can develop
$A_{[1,2,3,i], \cdot}$  by the first or the second row with coeffcients 
$0, s_i,0 $

By the  above remarks, we can assume that the minumum in $S_i$ is 
attained  twice by two entries in the same row.


Observe that if the minumum $s_i$ in $S_i$ is 
attained  twice by two entries in the same row and this row is not the last 
one, then we can easily conclude, because, if for instance this row is 
$A_{(2)}$,  we can develop $A_{[1,2,3,i],\cdot}$
by $A_{(2)}$  with coefficients of order $0,s_i,0$.

So we can assume that, for every $ i=5,...,g$,   the minimum
 $s_i$ in $S_i$ is attained in the last row.

So 
 $A= \left(\begin{array}{ccccc}
0 & 0 & c & c' & y \\
0 & 0 & d & d' & 0 \\
a & a' & 0 & 0 & 0 \\
b & b' & 0 & 0 & x \\
\cdot & \cdot & s_5 & s_5 & 0 \\
\vdots & \vdots & \vdots & \vdots & \vdots \\
\cdot & \cdot & s_g & s_g & 0 
\end{array}
\right)$

If $ c=c'$ and $d=d'$, we can conclude at once.

If $c=c' $ and $ d \neq d'$, we can suppose for symmetry that $d>d'$ and we 
can develop $ A_{\cdot , [2,3,4,5]}$ by $ A^{(4)}$ with coefficients 
$ y+d', 0 , d' $. 

If $d=d' $ and $ c \neq c'$, we can suppose for symmetry that $c>c'$ and we 
can develop $ A_{\cdot , [2,3,4]}$ by $ A^{(3)}$ with coefficients 
$ c', 0 $ if $ c' >d$, while we can develop  $ A_{\cdot , [2,3,4,5]}$ by 
$ A^{(4)}$ with coefficients 
$ c', 0, c' $ if $ c'  \leq d$. 

So we can suppose $d \neq d' $ and $ c \neq c'$

Suppose  $min \{d,d',c,c'\} \in \{d,d'\}$. For symmetry we can suppose
it is $d$. If $c=d $ or $c'=d$, 
we can develop  $ A_{\cdot, [2,3,4]}$ 
by $ A^{(3)}$ 
column with coeffficients of order respectively $ d ,0$.
 Thus we can suppose $c>d $ and $c'>d $;
we can develop  $ A_{\cdot, [2,3,4,5]}$ 
by $ A^{(3)}$ 
column with coeffficients of order respectively $ min \{c,c',d+y\} ,0,d$.

Analogously if $min \{d,d',c,c'\} \in \{c,c'\}$.

\underline{Case b)}
We can suppose $$A= \left(\begin{array}{ccccc}
0 & 0 & 0 & \ast & \ast \\
0 & 0 & \ast & 0 & \ast \\
. & . & 0 & \ast & 0 \\
. & . & \ast & 0 & 0 \\
\vdots & \vdots & \vdots & \vdots & \vdots
\end{array}
\right)$$
if not, u.s.r.c., we are in case a).
Let $R= A_{[34], [12]}$.  

If all entries of  $R$ are nonzero,
develop $A_{[1,2,3,4],\cdot } $ by the first row 
with coefficients of order $0,0,0$ (apply Remark \ref{remquadsbar}).
 
If in $R$ there are two zeroes in the same row we can conclude by Lemma
\ref{4linedup}.

If in $R$ there are exactly two zeroes but not in the same row or column, 
we can suppose that $ R= \left(\begin{array}{cc}
0 & \ast  \\
\ast & 0
\end{array}
\right)$.
Observe that in this case, for any $j=1,...,5$, the matrix $ A_{[1234], 
\hat{j}}$ has rank $3$ and in every of its rows there are exactly two 
zeroes, thus for any $ i=5,..., g$, for any $j=1,..., 5$, in $A_{(i), 
\hat{j}}  $ there are at least two zeroes; 
thus in   $A_{(i)}$ there are at least $3$ zeroes; we can 
suppose exactly $3$ by Lemma \ref{4linedup}. So we can conclude by Lemma
\ref{caso3322} (distinguish into the  cases: at least one between
 $A_{i,1}$ and $A_{i,3}$ is $0$, both are nonzero).

If in $R$ there are exactly two zeroes in the same column, as above,
 we can suppose that for any $ i=5,..., g$, 
 in   $A_{(i)}$ there are exactly $3$ zeroes.
We can conclude by Lemma \ref{caso3322}, except if $ A_{(i)}= 
(0,\ast , 0, \ast ,0) $ or $ A_{(i)}= 
(0, 0, \ast, \ast ,0) $. But if, for any $ i=5,..., g$, $ A_{(i)}$ is like 
this, then $ A^{(1)}=0$, and then we can develp $A$ by the first column.

Suppose in $R$ there is only one zero; we can suppose that it is in the
 first column; we can observe that 
 for any $j=1,...,3$, the matrix $ A_{[1234], 
\hat{j}}$ has rank $3$ and in every of its rows there are exactly two 
zeroes, thus for any $ i=5,..., g$, for any $j=1,..., 3$, in $A_{(i), 
\hat{j}}  $ there are at least two zeroes. We can conclude in any case 
by Lemma \ref{caso3322}, except if $A_{(i) }= (0,\ast , 0,0, \ast ) $ 
or $ A_{(i)}=  (0, 0, \ast, \ast ,0) $; in these two cases we can conclude 
by developping $A$ by the first row with coefficients of order $0$ by 
Remark \ref{sbis}.

\underline{Case c)}  
We can suppose  $$A= \left(\begin{array}{ccccc}
0 & 0 & 0 & \ast & \ast \\
0 & 0 & \ast & 0 & \ast \\
. & . & 0 & 0 & 0 \\
. & . & \ast & \ast & 0 \\
\vdots & \vdots & \vdots & \vdots & \vdots
\end{array}
\right)$$
if not, we are in case a).
By Lemma \ref{4linedup} we can suppose also 
$$A= \left(\begin{array}{ccccc}
0 & 0 & 0 & \ast & \ast \\
0 & 0 & \ast & 0 & \ast \\
\ast & \ast & 0 & 0 & 0 \\
. & . & \ast & \ast & 0 \\
\vdots & \vdots & \vdots & \vdots & \vdots
\end{array}
\right)$$

If $$A= \left(\begin{array}{ccccc}
0 & 0 & 0 & \ast & \ast \\
0 & 0 & \ast & 0 & \ast \\
\ast & \ast & 0 & 0 & 0 \\
\ast & \ast & \ast & \ast & 0 \\
\vdots & \vdots & \vdots & \vdots & \vdots
\end{array}
\right)$$
 we can develop 
$A_{[1234], \cdot}$ by the first row with coefficients $0,0,0$
 (apply Remark \ref{remquadsbar}).

So we can suppose $$A= \left(\begin{array}{ccccc}
0 & 0 & 0 & \ast & \ast \\
0 & 0 & \ast & 0 & \ast \\
\ast & \ast & 0 & 0 & 0 \\
. & 0 & \ast & \ast & 0 \\
\vdots & \vdots & \vdots & \vdots & \vdots
\end{array}
\right)$$
We can conclude by Lemma \ref{caso3322}.
\hfill \framebox(7,7)

\bigskip

{\bf Note.} In the contemporary joint paper  \cite{CJR}
a much more elegant proof of Theorem \ref{mio} is exhibited.

We also notice that the first version of this eprint (arXiv:0712.3007v1, 
December 2007) was wrong.

\end{document}